\documentclass[12pt]{amsart}
\usepackage{amsmath, amssymb, latexsym, amsthm, bm}
\usepackage{mathrsfs}
\usepackage[all]{xy}

\newcommand{\Pa}[9]{\bibitem{#1} {#2}, \emph{#3}, {#4} \textbf{#5} ({#6}), {#7}--{#8}.}

\newcommand{\ed}{

\end{document}
}

\newcommand{\arx}[1]{\texttt{arxiv.org/math/#1}}
\newcommand{\bq}{\begin{quote}}
\newcommand{\eq}{\end{quote}}

\newcommand{\inv}{^{-1}}
\newcommand{\Cantor}{{\{0,1\}^\N}}

\newcommand{\N}{\mathbb{N}}
\newcommand{\NN}{{\N^{\N}}}

\newcommand{\PN}{{P(\N)}}
\newcommand{\roth}{{[\N]^{\aleph_0}}}

\newcommand{\seq}[1]{\{#1\}_{n\in\N}}
\newcommand{\sseq}[1]{\setseq{#1}}
\newcommand{\setseq}[1]{\{#1 : n\in\N\}}

\newcommand{\cJ}{\mathcal{J}}
\newcommand{\scrA}{\mathscr{A}}
\newcommand{\scrB}{\mathscr{B}}

\newcommand{\B}{\mathcal{B}}

\newcommand{\BG}{\B_\Gamma}

\newcommand{\CO}{C_\Omega}

\newcommand{\cF}{\mathcal{F}}

\newcommand{\cO}{\mathcal{O}}

\newcommand{\R}{\mathbb{R}}
\newcommand{\cU}{\mathcal{U}}
\newcommand{\Union}{\bigcup}
\newcommand{\cV}{\mathcal{V}}
\newcommand{\cW}{\mathcal{W}}

\newcommand{\Impl}{\Rightarrow}
\long\def\forget#1\forgotten{}

\newcommand{\fd}{\mathfrak{d}}

\newcommand{\oo}{\infty}

\newcommand{\fp}{\mathfrak{p}}

\newcommand{\w}{\omega}

\newcommand{\nin}{\not\in}

\newcommand{\sbst}{\subseteq}

\newcommand{\sm}{\setminus}
\newcommand{\as}{\subseteq^*}
\newcommand{\rest}{\restriction}

\newcommand{\non}{\mathsf{non}}

\newtheorem{thm}{Theorem}
\newtheorem{prop}[thm]{Proposition}

\newtheorem{prob}[thm]{Problem}
\newtheorem{lem}[thm]{Lemma}
\newtheorem{cor}[thm]{Corollary}
\newtheorem{conj}[thm]{Conjecture}
\theoremstyle{definition}
\newtheorem{defn}[thm]{Definition}
\theoremstyle{remark}

\newtheorem{exam}[thm]{Example}
\newcommand{\be}{\begin{enumerate}}
\newcommand{\ee}{\end{enumerate}}
\newcommand{\bi}{\begin{itemize}}
\newcommand{\itm}{\item}
\newcommand{\ei}{\end{itemize}}

\forget
\forgotten

\newcommand{\sone}{\mathsf{S}_1}
\newcommand{\sfin}{\mathsf{S}_{fin}}
\newcommand{\ufin}{\mathsf{U}_{fin}}

\newcommand{\capinf}{\mbox{$\bigcap\nolimits_{\oo}$}}
\newcommand{\strongbinom}[2]{\strongbinomlight{#1}{\bm{#2}}}
\newcommand{\strongbinomlight}[2]{\binom{#1}{{#2}_\oo}}

\author{Boaz Tsaban}
\thanks{Supported by the Koshland Center for Basic Research.}
\address{Boaz Tsaban, Department of Mathematics,
Weizmann Institute of Science,
Rehovot 76100,
Israel}
\email{boaz.tsaban@weizmann.ac.il}
\urladdr{http://www.cs.biu.ac.il/\~{}tsaban}

\title{A new selection principle}

\begin{document}

\begin{abstract}
Motivated by a recent result of Sakai, we define a new selection
operator for covers of topological spaces, inducing new selection
hypotheses, and initiate a systematic study of the new
hypotheses.
Some intriguing problems remain open.
\end{abstract}

\maketitle

\section{Subcovers with strong covering properties}

We say that $\cU$ is a \emph{cover} of a set $X$ if $X\nin\cU$
and $X=\Union\cU$.

\begin{defn}\label{Ainf}
For a family $\scrA$ of covers of a set $X$, $\scrA_\oo$ is the family of all
$\cU$ such that there exist infinite sets $\cU_n\sbst\cU$, $n\in\N$, with
$\setseq{\bigcap\cU_n}\in\scrA$.
\end{defn}

For topological spaces $X$, various special families of covers have been
extensively studied in the literature, in a framework called \emph{selection
principles}, see the surveys \cite{LecceSurvey, KocSurv, ict}.
The main types of covers are defined as follows.
Let $\cU$ be a cover of $X$. $\cU$ is an \emph{$\omega$-cover} of $X$ if
each finite $F\sbst X$ is contained in some $U\in\cU$.
$\cU$ is a \emph{$\gamma$-cover} of $X$ if
$\cU$ is infinite, and each $x\in X$ belongs to all but finitely many $U\in\cU$.

Let the boldfaced symbols $\bm{\cO}$, $\bm{\Omega}$, $\bm{\Gamma}$
denote the families of all covers, $\omega$-covers, and $\gamma$-covers, respectively.
Then
$$\bm{\Gamma}\sbst \bm{\Omega}\sbst \bm{\cO}.$$
Also, let $\cO$, $\Omega$, $\Gamma$ denote the corresponding families of \emph{open} covers.

For a space $X$ and collections $\scrA,\scrB$ of covers of $X$,
the following property may or may not hold:
\bi
\itm[$\binom{\scrA}{\scrB}$:] Every member of $\scrA$ has a subset which is a member of $\scrB$.
\ei

Sakai \cite{Sakai03} proved that for a Tychonoff space $X$,
a local property called the \emph{Pytkeev property} holds in the function space $C_p(X)$
if, and only if, $X$ satisfies $\strongbinom{\check{\Omega}}{\Omega}$, where
$\check{\Omega}$ is a certain subclass of $\Omega$. It is open whether
$\check{\Omega}$ can be replaced by $\Omega$ \cite{Sakai03}.
Motivated by this, Simon and the present author proved that for Lindel\"of
spaces $X$,
$C_p(X)$ satisfies the Pytkeev property if, and only if,
$X$ is zero-dimensional and satisfies $\strongbinom{\CO}{\Omega}$,
where $\CO$ is the collection of all \emph{clopen} $\omega$-covers of $X$  \cite{Pyt}.
This motivates the study of additional properties
of the form $\strongbinomlight{\scrA}{\scrB}$.

\begin{defn}
A family $\scrB$ of open covers of $X$ is \emph{surjectively derefinable}
if for each $\cU\in\scrB$ and each $f:\cU\to P(X)\sm\{X\}$
such that for each $U\in\cU$ $f(U)$ is open and contains $U$,
$\{f(U) : U\in\cU\}\in\scrB$.
A similar definition applies to families of Borel covers, clopen covers,
etc.
\end{defn}

\begin{exam}
$\bm\cO$, $\bm\Omega$, and $\bm\Gamma$ are surjectively derefinable.
For the latter we must explain why a (surjective) derefinement of a $\gamma$-cover
is infinite, and this follows from the fact that it is an $\omega$-cover.
\end{exam}

\begin{lem}\label{surde}
Assume that $\scrB$ is a surjectively derefinable family of covers of $X$.
Then
$$\strongbinomlight{\scrA}{\scrB}\Impl \binom{\scrA}{\scrB}.$$
\end{lem}
\begin{proof}
Assume that $\cU\in\scrA$. By the assumption, there are
infinite $\cU_1,\cU_2,\dots\sbst\cU$ such that $\cV=\setseq{\bigcap\cU_n}\in\scrB$.
For each $n$, choose $f(\bigcap\cU_n)\in\cU_n$.
As $\scrB$ is surjectively derefinable, $\cW=\setseq{f(\bigcap\cU_n)}\in\scrB$.
Clearly, $\cW\sbst\cU$.
\end{proof}

The converse need not hold. For example, $\binom{\CO}{\bm\Omega}$
always holds, whereas $\strongbinom{\CO}{\Omega}$ need not,
as explained above. More examples will follow in the sequel.

\begin{prop}\label{gaga}
Every space satisfies $\strongbinom{\Gamma}{\Gamma}$.
\end{prop}
\begin{proof}
Assume that $\cU\in\Gamma$.
We may assume that $\cU$ is countable (since an infinite subset of a $\gamma$-cover is again a
$\gamma$-cover). Enumerate $\cU=\setseq{U_n}$ bijectively,
and take $\cU_n=\{U_k : k\ge n\}$ for each $n$.
Then $\setseq{\bigcap\cU_n}\in\Gamma$.
\end{proof}

$\binom{\Omega}{\bm\Gamma}$ is the classical $\gamma$-property \cite{GN}.

\begin{cor}\label{omgam}
$\binom{\Omega}{\bm\Gamma}=\strongbinom{\Omega}{\Gamma}$.
\end{cor}
\begin{proof}
By Proposition \ref{gaga} and Lemma \ref{surde},
$$\binom{\Omega}{\bm\Gamma}=\binom{\Omega}{\bm\Gamma}\cap\strongbinom{\Gamma}{\Gamma}
 = \strongbinom{\Omega}{\Gamma},$$
 the last equation being self-evident.
\end{proof}

What about the other properties? $\strongbinom{\cO}{\cO}$ never
holds, since any $T_1$ space with more than $1$ element has a
finite open cover. Taking the above results into account, only
$\strongbinom{\Omega}{\cO}$ and $\strongbinom{\Omega}{\Omega}$ are
potentially new. It turns out that even the formally weaker
property $\strongbinom{\Omega}{\cO}$ is quite restrictive.
According to Borel, a set $X\sbst\R$ has \emph{strong measure
zero} if for each sequence of positive reals $\seq{\epsilon_n}$,
there exists a cover $\seq{I_n}$ of $X$ such that for each $n$,
the diameter of $I_n$ is smaller than $\epsilon_n$. It was
established by Laver that consistently, all strong measure zero
sets of reals are countable. The following theorem of Miller
is essentially proved in \cite{Pyt}.

\begin{thm}[Miller]\label{SMZ}
If $X\sbst\R$ and $X$ satisfies $\strongbinom{\CO}{\cO}$,
then $X$ has strong measure zero.
\end{thm}
\begin{proof}
By standard arguments \cite{prods}, we may assume that $X\sbst\{0,1\}^\N$.
It suffices to prove that for each increasing sequence $\seq{k_n}$ of natural
numbers, there are for each $n$ elements $s^n_m\in\{0,1\}^{k_n}$, $m\le n$, such
that $X=\Union_n ([s^n_1]\cup\dots\cup [s^n_n])$. (One can allow $n$ sets of diameter $\epsilon_n$
in the original definition of strong measure zero by moving to an appropriate
subsequence of the original sequence $\seq{\epsilon_n}$.)

For each $n$, let
$$\cU_n = \{[s_1]\cup\dots\cup [s_n] : s_1,\dots,s_n\in\{0,1\}^{k_n}\},$$
and take $\cU=\Union_n\cU_n$.
$\cU$ is a clopen $\omega$-cover of $X$. By $\strongbinom{\CO}{\cO}$,
there are infinite subsets $\cV_1,\cV_2,\dots$ of $\cU$, such that $\setseq{\bigcap\cV_n}$ is
a cover of $X$.
As each $\cV_n$ is infinite and each $\cU_n$ is finite, we can find $m_1$ and $V_1\in\cV_1\cap\cU_{m_1}$,
$m_2>m_1$ and $V_2\in\cV_2\cap\cU_{m_2}$, etc. Then $\{V_n : n\in\N\}$ is
a cover of $X$, and the sets $V_n$ are as required in the first paragraph
of this proof.
\end{proof}

However, we have the following.

\begin{conj}\label{conj}
The Continuum Hypothesis implies:
\be
\itm There is a set of reals $X$ satisfying
$\strongbinom{\Omega}{\Omega}$ but not $\binom{\Omega}{\bm\Gamma}$; and
\itm There is a set of reals $X$ satisfying
$\strongbinom{\Omega}{\cO}$ but not $\strongbinom{\Omega}{\Omega}$.
\ee
\end{conj}

Conjecture \ref{conj}(1) implies, if true, a negative answer to Sakai's Question 4.7 in \cite{Sakai06}.
We will show that critical cardinalities (defined below) cannot be used to prove
the consistency of items (1) and (2) of Conjecture \ref{conj}.

$X$ is an \emph{$\Omega$-Lindel\"of} space if
each open $\omega$-cover of $X$ contains a countable
$\omega$-cover of $X$. For Tychonoff spaces this is equivalent to:
All finite powers of $X$ are Lindel\"of \cite{GN}.
Separable zero-dimensional
metrizable spaces are homeomorphic to subsets of $\R$, and are
thus $\Omega$-Lindel\"of.
Recall that a family $\cF\sbst\roth$ is \emph{centered}
if the intersection of each finite subset of $\cF$ is infinite.
$\cF$ is \emph{free} if $\bigcap\cF=\emptyset$.
$A\sbst\N$ is a \emph{pseudo-intersection} of $\cF$ if
$A$ is infinite and for each $B\in\cF$, $A\as B$ (that is,
$A\sm B$ is finite). $\roth$ inherits its topology from
$P(\N)$, whose topology is defined by identifying $P(\N)$ with
$\Cantor$.

\begin{thm}\label{OmOchar}
For $\Omega$-Lindel\"of spaces $X$, the following are equivalent:
\be
\itm $X$ satisfies $\strongbinom{\CO}{\cO}$;
\itm For each continuous free centered image $\cF$ of $X$ in $\roth$,
$\cF=\Union_n\cF_n$ where each $\cF_n$ has a pseudo-intersection.
\ee
\end{thm}
\begin{proof}
$(1\Impl 2)$ Assume that $\Psi:X\to\roth$ is continuous and that
its image $\cF$ is free and centered. For each $n$, let $U_n = \{x : n\in\Psi(x)\}$.
$\cU=\sseq{U_n}$ is a clopen $\omega$-cover of $X$. Choose infinite
$\cU_n\sbst\cU$, $n\in\N$, such that $\sseq{\bigcap\cU_n}$ is a cover of $X$,
and set $A_n=\{m : U_m\in\cU_n\}$, and $\cF_n = \{I\in\cF : A_n\sbst I\}$.
For each $I\in\cF$, let $x\in X$ be such that $I=\Psi(x)$.
Choose $n$ such that $x\in \bigcap\cU_n$. Then for each $m\in A_n$,
$x\in U_m$ and therefore $m\in\Psi(x)=I$, that is, $I\in\cF_n$.

$(2\Impl 1)$ Assume that $\cU$ is a clopen $\omega$-cover of $X$.
Since $X$ is $\Omega$-Lindel\"of, we may assume that $\cU$ is countable.
Fix a bijective enumeration $\cU=\sseq{U_n}$. As the sets $U_n$ are clopen,
the Marczewski function $\mu:X\to \PN$ defined by
$\mu(x)=\{n : x\in U_n\}$ is continuous.
Since $\cU$ is an $\omega$-cover of $X$, the image $\cF$ of $\mu$
is a free centered subset of $\roth$ \cite{ict}.
Let $\cF=\Union_n\cF_n$ be as in (2).
For each $n$, let $A_n$ be a pseudo-intersection of $\cF_n$.
Take $\cU_{n,m} = \{U_k : m\le k\in A_n\}\sbst\cU$.
Then $\{\bigcap \cU_{n,m} : m,n\in\N\}$ is a cover of $X$.
\end{proof}

The minimal cardinality of a centered $\cF\sbst\roth$
such that there is \emph{no} partition $\cF=\Union_n\cF_n$ where each
$\cF_n$ has a pseudo-intersection is equal to $\fp$ \cite{Pyt}.

The \emph{critical cardinality} of a nontrivial family $\cJ$ of sets of reals
is
$$\non(\cJ)=\min\{|X| : X\sbst\R\mbox{ and }X\nin\cJ\}.$$
\begin{cor}\label{p}
$\non(\strongbinom{\Omega}{\Omega})=\non(\strongbinom{\Omega}{\cO})=\fp$.
\end{cor}
\begin{proof}
$\non(\binom{\Omega}{\bm\Gamma})=\fp$ \cite{GM}, and by the implications among the properties,
$\fp\le\non(\strongbinom{\Omega}{\Omega})\le\non(\strongbinom{\Omega}{\cO})\le\non(\strongbinom{\CO}{\cO})$.
By Theorem \ref{OmOchar} and the above-mentioned result of \cite{Pyt},
$\non(\strongbinom{\CO}{\cO})\le\fp$.
\end{proof}

\begin{prop}\label{powers}
If all finite powers of $X$ satisfy $\strongbinom{\Omega}{\cO}$, then
$X$ satisfies $\strongbinom{\Omega}{\Omega}$.
\end{prop}
\begin{proof}
If $\cU$ is an open $\w$-cover of $X$, then for each $k$,
$\cU^k:=\{U^k : U\in\cU\}$ is an open $\w$-cover of $X^k$.
Take infinite $\cV_{k,n}\sbst\cU$ such that $\sseq{\bigcap\cV_{k,n}^k}$ is a cover of $X^k$.
Then each $k$-element subset of $X$ is contained in some member of
$\sseq{\bigcap\cV_{k,n}}$, and therefore $\{\bigcap\cV_{k,n} : n,k\in\N\}$
is an $\w$-cover of $X$.
\end{proof}

A subtle technical problem prevents us from using the methods
of \cite{coc2} to obtain the converse implication.

\begin{prob}
Is the converse implication in Proposition \ref{powers} provable?
\end{prob}

Additional results concerning $\strongbinom{\Omega}{\Omega}$ can be found in \cite{Sakai06}.

\section{A new selection principle}

Fix a topological space $X$, and let
$\scrA$ and $\scrB$ each be a collection of covers of $X$.
The following selection principles, which $X$ may or may not satisfy,
were introduced in \cite{coc1} to generalize a variety of classical
properties, and were extensively studied in the literature
(see the surveys \cite{LecceSurvey, KocSurv, ict}).

\begin{itemize}
\item[$\sone(\scrA,\scrB)$:]
For each sequence $\seq{\cU_n}$ of members of $\scrA$,
there exist members $U_n\in\cU_n$, $n\in\N$, such that $\setseq{U_n}\in\scrB$.
\item[$\sfin(\scrA,\scrB)$:]
For each sequence $\seq{\cU_n}$
of members of $\scrA$, there exist finite
subsets $\cF_n\sbst\cU_n$, $n\in\N$, such that $\Union_{n\in\N}\cF_n\in\scrB$.
\item[$\ufin(\scrA,\scrB)$:]
For each sequence $\seq{\cU_n}$ of members of $\scrA$
which do not contain a finite subcover,
there exist finite subsets $\cF_n\sbst\cU_n$, $n\in\N$,
such that $\setseq{\cup\cF_n}\in\scrB$.
\end{itemize}

We introduce the following new selection principle,
which is a selective version of $\strongbinomlight{\scrA}{\scrB}$:
\bi
\itm[$\capinf(\scrA,\scrB):$] For each sequence $\seq{\cU_n}$ of elements of $\scrA$,
there is for each $n$ an infinite set $\cV_n\sbst\cU_n$, such that
$\setseq{\bigcap\cV_n}\in\scrB$.
\ei
Note that if $\scrA$ contains a finite element, then $\capinf(\scrA,\scrB)$
automatically fails.

Since the sequence $\seq{\cU_n}$ is allowed to be constant,
the following holds.
\begin{lem}
$\capinf(\scrA,\scrB)\Impl \strongbinomlight{\scrA}{\scrB}$.\hfill\qed
\end{lem}

The following is easy to verify.

\begin{prop}\label{capinfbasic}
Assume that $\scrB$ is a surjectively derefinable family of covers of $X$.
Then
$$\capinf(\scrA,\scrB)\Impl \sone(\scrA,\scrB).\qed$$
\end{prop}

As $\capinf(\Gamma,\bm\Gamma)\Impl\sone(\Gamma,\bm\Gamma)$, and $\sone(\Gamma,\bm\Gamma)$ is rather restrictive
(e.g., every set of reals satisfying it is perfectly meager),
it follows that $\capinf(\Gamma,\bm\Gamma)$ is strictly stronger than $\strongbinom{\Gamma}{\Gamma}$.

\begin{thm}\label{capinfgg}
$\capinf(\Gamma,\bm\Gamma)=\sone(\Gamma,\bm\Gamma)$.
\end{thm}
\begin{proof}
Assume that $X$ satisfies $\sone(\Gamma,\bm\Gamma)$.
We will prove that $X$ satisfies $\capinf(\Gamma,\bm\Gamma)$.
The trick we use comes from the context of local properties,
in which we learned it from Scheepers.

Assume that $\cU_n$, $n\in\N$, are open $\gamma$-covers of $X$. We
may assume that they are all countable and that the sets $\cU_n$, $n\in\N$, are
pairwise disjoint.

Fix a surjection
$f:\N\to\N$ such that for each $n$, $f\inv(n)$ is infinite. For a
countable bijectively enumerated set $\cF=\sseq{U_n}$ and $m\in\N$, define
$\cF(m)=\{U_n : n\ge m\}$. Fix a bijective enumeration for each of the
covers $\cU_n$, and apply $\sone(\Gamma,\bm\Gamma)$ to the
sequence $\cU_{f(n)}(n)$, $n\in\N$, to obtain sets
$U_n\in\cU_{f(n)}(n)$ such that $\sseq{U_n}$ is a $\gamma$-cover
of $X$.

For each $n$, take $\cV_n = \{U_m : f(m)=n\}\sbst\cU_n$.
Each $U_m\in\cV_n$ can belong to only
finitely many $\cU_{f(k)}(k)$ with $f(k)=n$, and cannot belong to
any $\cU_{f(k)}(k)$ with $f(k)\neq n$ (because $\cU_n\cap\cU_{f(k)}=\emptyset$).
In particular, $\cV_n$ is infinite.

$\cV=\sseq{\bigcap\cV_n}$ is a $\gamma$-cover of $X$:
For each $x\in X$, $x\in U_m$ for all large enough $m$,
and as $n\to\infty$, $\min f\inv(n)\to\infty$ either.
This also shows that $\cV$ is an $\omega$-cover of $X$, and
thus $\cV$ is infinite.
\end{proof}

$\sone(\Omega,\bm\Gamma)=\binom{\Omega}{\bm\Gamma}$ \cite{GN}.

\begin{cor}\label{gamma2}
$\capinf(\Omega,\bm\Gamma)=\sone(\Omega,\bm\Gamma)$.
\end{cor}
\begin{proof}
By Theorem \ref{capinfgg} and easy reasoning,
$$\capinf(\Omega,\bm\Gamma)=\binom{\Omega}{\bm\Gamma}\cap\capinf(\Gamma,\bm\Gamma)= \binom{\Omega}{\bm\Gamma}\cap
\sone(\Gamma,\bm\Gamma)=\sone(\Omega,\bm\Gamma).\qedhere$$
\end{proof}

Exactly the properties in Figure \ref{capsurv} remain to be explored.

\begin{figure}[!htp]
$\xymatrix{
\capinf(\Gamma,\bm\Omega)\ar[r] & \capinf(\Gamma,\bm\cO)\\
\capinf(\Omega,\bm\Omega)\ar[u]\ar[r] & \capinf(\Omega,\bm\cO)\ar[u]
}$
\caption{The surviving properties}\label{capsurv}
\end{figure}

\begin{prob}\label{classi}
Is any of the properties in Figure \ref{capsurv} equivalent to a
classical selection hypothesis?
\end{prob}

\begin{prob}
Can any implication be added to Figure \ref{capsurv}?
\end{prob}

For $\capinf(\Omega,\bm\Omega)$ and $\capinf(\Omega,\bm\cO)$,
Problem \ref{classi} is closely related to Conjecture \ref{conj},
because these properties are sandwiched between $\binom{\Omega}{\bm\Gamma}$
and $\strongbinom{\Omega}{\cO}$.

For the remaining two properties, we have a partial answer for Problem \ref{classi}.
Let $\BG$ denote the family of \emph{countable Borel} $\gamma$-covers of $X$.

\begin{thm}\label{borelcap}
~\be
\itm $\sone(\BG,\bm\cO)=\capinf(\BG,\bm\cO)$.
\itm $\sone(\BG,\bm\Omega)=\capinf(\BG,\bm\Omega)$.
\itm $\sone(\BG,\bm\cO)\Impl \capinf(\Gamma,\bm\cO)\Impl \sone(\Gamma,\bm\cO)$.
\itm $\sone(\BG,\bm\Omega)\Impl \capinf(\Gamma,\bm\Omega)\Impl \sone(\Gamma,\bm\Omega)$.
\ee
\end{thm}
\begin{proof}
We only prove the implications which do not follow from
Proposition \ref{capinfbasic}.

(1) Assume that $X$ satisfies $\sone(\BG,\bm\cO)$ and
$\cU_n$, $n\in\N$, are countable Borel $\gamma$-covers of $X$.
Enumerate bijectively, for each $n$, $\cU_n=\{U^n_m : m\in\N\}$.
Define $\Psi:X\to\NN$ by
$$\Psi(x)(n)=\min\{m : (\forall k\ge m)\ x\in U^n_k\}.$$
Since $\Psi$ is Borel and $X$ satisfies $\sone(\BG,\bm\cO)$, $\Psi[X]$ is not dominating \cite{CBC}.
Let $g\in\NN$ be a witness for that.
Take $\cV_n = \{U^n_m : m\ge g(n)\}$. For each $x\in X$, there are infinitely many $n$ such that
$\Psi(x)(n)\le g(n)$, and therefore $x\in\bigcap \cV_n$.

(2) is similar, here $\Psi[X]$ is not finitely dominating \cite{CBC}, and this is what we need.

(3) and (4) follow from (1) and (2), respectively, because $\gamma$-covers may be assumed
to be countable.
\end{proof}

\begin{cor}
~\be
\itm $\non(\capinf(\Omega,\bm\Omega))=\non(\capinf(\Omega,\bm\cO))=\fp$.
\itm $\non(\capinf(\Gamma,\bm\Omega))=\non(\capinf(\Gamma,\bm\cO))=\fd$.
\ee
\end{cor}
\begin{proof}
(1) follows from Corollaries \ref{gamma2} and \ref{p} and the implications among the properties,
together with $\non(\sone(\Omega,\bm\Gamma))=\fp$ \cite{GM}.

(2) follows from Theorem \ref{borelcap} and the fact that the critical cardinalities of the
Borel version of the classical principles are the same as in their open version, and are both $\fd$
\cite{coc2, CBC}.
\end{proof}

There are some additional interesting connections between the new and the classical
selection principles.

\begin{thm}[Sakai \cite{Sakai06}]\label{sak}
$\sfin(\Omega,\bm\Omega)\cap\strongbinom{\Omega}{\Omega}\Impl\sone(\Omega,\bm\Omega)$.
\end{thm}
\begin{proof}
For completeness, we give a proof.

Clearly,
$$\sfin(\Omega,\bm\Omega)\cap\strongbinom{\Omega}{\Omega}=\sfin(\Omega,\bm\Omega_\oo).$$
It therefore suffices to show that $\sfin(\Omega,\bm\Omega_\oo)$ implies $\sone(\Omega,\bm\Omega)$.
Indeed, assume that $\cU_n$, $n\in\N$, are open $\omega$-covers of $X$.
Choose finite $\cF_n\sbst\cU_n$, $n\in\N$, such that $\Union_n\cF_n\in\bm\Omega_\oo$ for $X$.

Take infinite $\cV_n\sbst\Union_n\cF_n$, $n\in\N$, such that
$\sseq{\bigcap\cV_n}$ is an $\omega$-cover of $X$.
To each $n$, assign $m_n$ such that $m_n$ is increasing with $n$ and $\cV_n\cap\cF_{m_n}\neq\emptyset$,
and choose any $U_{m_n}\in\cV_n\cap\cF_{m_n}$.
For $k\nin\sseq{m_n}$ choose any $U_k\in\cU_k$.
As $\sseq{\bigcap\cV_n}$ refines $\sseq{U_n}$, $\sseq{U_n}$ is an $\omega$-cover of $X$.
\end{proof}

\begin{cor}
$\sfin(\Omega,\bm\Omega_\oo)=\sone(\Omega,\bm\Omega_\oo)$.
\end{cor}
\begin{proof}
By Sakai's Theorem \ref{sak}:
\begin{eqnarray*}
\lefteqn{\sfin(\Omega,\bm\Omega_\oo) =}\\
& = & \sfin(\Omega,\bm\Omega)\cap\strongbinom{\Omega}{\Omega}=
\sone(\Omega,\bm\Omega)\cap\strongbinom{\Omega}{\Omega}
 = \sone(\Omega,\bm\Omega_\oo).\qedhere
\end{eqnarray*}
\end{proof}

The following result is inspired by results from \cite{PytII}.

\begin{thm}\label{basic}
For Lindel\"of zero-dimensional spaces, $\ufin(\cO,\Gamma)=\ufin\allowbreak(\cO,\bm\cO_\oo)$.
\end{thm}
\begin{proof}
Note that $\Gamma\sbst\bm\cO_\oo$. We therefore prove
that $\ufin(\cO,\allowbreak\bm\cO_\oo)$ implies $\ufin(\cO,\Gamma)$.
Assume that $X$ is Lindel\"of zero-dimensional, and satisfies $\ufin(\cO,\bm\cO_\oo)$. It suffices to prove that
every continuous image of $X$ in $\NN$ is bounded \cite{Rec94}.

Assume that $Y$ is a continuous image of $X$ in $\NN$.
We may assume that all elements of $Y$ are increasing functions.
If there is an infinite $I\sbst\N$ such that $\{f\rest I : f\in Y\}$ is
bounded, then $Y$ is bounded. We therefore assume that there is $k$
such that for each $n\ge k$, $\{f(n) : f\in Y\}$ is infinite.

For each $n\ge k$, let $\cU_n=\{U^n_m : m\in\N\}$, where
$U^n_m=\{f\in Y : f(n)\le m\}$ for each $m$.
$\cU_n$ does not contain $Y$ as an element.
Thus, there are finite sets $\cF_n\sbst\cU_n$, $n\ge k$, such that
$\cV=\sseq{\Union\cF_n}\in\bm\cO_\oo$. We may assume that each
$\cF_n$ is nonempty.
For each $n$, the sets $U^n_m$ are increasing
with $m$, and therefore there is $g(n)\in\N$ such that $\Union\cF_n=U^n_{g(n)}$.

Let $\cV_m$, $m\in\N$, be infinite subsets
of $\cV=\sseq{U^n_{g(n)}}$ such that $\{\bigcap\cV_m : m\in\N\}$ is a cover of $Y$.
For each $m$, let $I_m = \{n\ge k : U^n_{g(n)}\in\cV_m\}$.
$I_m$ is infinite, and $\{f\rest I_m : f\in \bigcap\cV_m\}$ is bounded by $g\rest I_m$. Thus,
$\bigcap\cV_m$ is bounded. It follows that $Y=\Union_m\bigcap\cV_m$ is bounded.
\end{proof}

Theorem \ref{basic} can be contrasted with the fact that $\ufin(\cO,\Gamma)\neq\ufin\allowbreak(\cO,\cO)$
\cite{SFH}.

\subsection*{Acknowledgments}
We thank Nadav Samet, Lyubomyr Zdomskyy, Dominic van der Zypen, and the referee for their useful
comments on this paper.
After we have submitted this paper for publication, we learned
of Ko\v{c}inac's earlier paper \cite{KocAlpha},
were a result very close to Theorem \ref{capinfgg} is proved using essentially
the same argument.
A complete classification of Ko\v{c}inac's properties in our context appears
in \cite{TsKocAlpha}.

\ed
\begin{thebibliography}{00}
\bibitem{GM}
F.\ Galvin and A.\ W.\ Miller,
\emph{$\gamma$-sets and other singular sets of real numbers},
Topology and its Applications \textbf{17} (1984),
145--155.

\bibitem{GN}
J.\ Gerlits and Zs.\ Nagy,
\emph{Some properties of $C(X)$, I},
Topology and its Applications \textbf{14} (1982),
151--161.

\bibitem{coc2}
W.\ Just, A.\ W.\ Miller, M.\ Scheepers, and P.\ J.\ Szeptycki,
\emph{The combinatorics of open covers II},
Topology and its Applications \textbf{73} (1996),
241--266.

\bibitem{KocSurv}
Lj.\ D.R.\ Ko\v{c}inac,
\emph{Selected results on selection principles},
in: \textbf{Proceedings of the 3rd Seminar on Geometry and Topology} (Sh.\ Rezapour, ed.),
July 15-17, Tabriz, Iran, 2004,
71--104.

\bibitem{KocAlpha}
Lj.\ D.R.\ Ko\v{c}inac,
\emph{Selection principles related to $\alpha_i$-properties},
Taiwanese Journal of Mathematics 12 (2008), 561--572.

\bibitem{Rec94}
I.\ Rec\l{}aw,
\emph{Every Luzin set is undetermined in the point-open game},
Fundamenta Mathematicae \textbf{144} (1994), 43--54.

\bibitem{Sakai03}
M.\ Sakai,
\emph{The Pytkeev property and the Reznichenko property in function spaces},
Note di Matematica \textbf{22} (2003),
43--52.

\bibitem{Sakai06}
M.\ Sakai,
\emph{Special subsets of reals characterizing local properties of function spaces},
in: \textbf{Selection Principles and Covering Properties in Topology} (L. Ko\v{c}inac, ed.),
Quaderni di Matematica \textbf{18},  Seconda Universita di Napoli, Caserta, 2006, 195--225.

\bibitem{coc1}
M.\ Scheepers,
\emph{Combinatorics of open covers I: Ramsey theory},
Topology and its Applications \textbf{69} (1996),
31--62.

\bibitem{LecceSurvey}
M.\ Scheepers,
\emph{Selection principles and covering properties in topology},
Note di Matematica \textbf{22} (2003),
3--41.

\bibitem{CBC}
M.\ Scheepers and B.\ Tsaban,
\emph{The combinatorics of Borel covers},
Topology and its Applications \textbf{121} (2002),
357--382.

\Pa{Pyt}{P. Simon and B. Tsaban}{On the Pytkeev property in spaces of continuous functions}{Proceedings of the American Mathematical Society}{136}{2008}{1125}{1135}{\arx{math}{GN}{0606270}}

\bibitem{ict}
B.\ Tsaban,
\emph{Some new directions in infinite-combinatorial topology},
in: \textbf{Set Theory} (J.\ Bagaria and S.\ Todor\v{c}evic, eds.),
Trends in Mathematics, Birkhauser, 2006, 225--255.

\bibitem{TsKocAlpha}
B.\ Tsaban,
\emph{On the Ko\v{c}inac $\alpha_i$ properties},
Topology and its Applications \textbf{155} (2007), 141--145.

\bibitem{prods}
B.\ Tsaban and T.\ Weiss,
\emph{Products of special sets of real numbers},
Real Analysis Exchange \textbf{30} (2004/5), 819--836.

\Pa{SFH}{B. Tsaban and L. Zdomskyy}{Scales, fields, and a problem of Hurewicz}{Journal of the European Mathematical Society}{10}{2008}{837}{866}{\arx{math}{GN}{0507043}}

\bibitem{PytII}
B.\ Tsaban and L.\ Zdomskyy,
\emph{On the Pytkeev property in spaces of continuous functions (II)},
Houston Journal of Mathematics \textbf{35} (2009), 563--571.

\end{thebibliography}
